\documentclass{article}
\usepackage{mlspconf}

\usepackage[vietnamese]{babel}

\usepackage{cite}
\usepackage{amsmath,amssymb,amsfonts}
\usepackage{graphicx}
\usepackage{textcomp}
\usepackage{xcolor}
\usepackage{comment}
\usepackage{url}

\usepackage[ruled,lined]{algorithm2e}
\usepackage{algpseudocode}

\usepackage{theorem}
\PassOptionsToPackage{normalem}{ulem}
\usepackage{ulem}

\newcommand{\bregmanloss}[1][\psi]{B_{#1}}
\newcommand{\R}{\mathbb{R}}
\newcommand{\N}{\mathbb{N}}
\newcommand{\eS}{\mathbb{S}}
\newcommand{\prox}[1]{\ensuremath{{\rm prox}_{#1}}}

\newcommand{\minimize}[1]{\ensuremath{\underset{#1}{{\rm minimize}}\,\,}}
\newcommand{\argmind}[1]{\ensuremath{\underset{#1}{{\rm argmin}}\,\,}}
\newcommand{\Idt}{\ensuremath{{\rm Id}}}

\newcommand{\ub}{\ensuremath{\boldsymbol{u}}}

\newcommand{\thetab}{\ensuremath{\boldsymbol{\theta}}}

\newcommand{\stkout}[1]{\ifmmode\text{\sout{\ensuremath{#1}}}\else\sout{#1}\fi}
\newcommand{\review}[1]{{\color{blue}#1}}

\newtheorem{theorem}{Theorem}[section]

\newtheorem{proposition}[theorem]{Proposition}

\copyrightnotice{979-8-3503-7225-0/24/\$31.00 \copyright 2024 IEEE}

\toappear{2024 IEEE International Workshop on Machine Learning for Signal Processing, Sept.\ 22--25, 2024, London, UK}

\title{A lifted Bregman strategy for training \\
unfolded proximal neural network Gaussian denoisers}
\name{%
    Xiaoyu Wang$^{\star}$%
    \qquad Martin Benning$^{\dagger}$%
    \qquad Audrey Repetti$^{\star}$%
}
\address{%
    $^{\star}$School of Mathematical
and Computer Sciences, Heriot-Watt University\\%
    $^{\dagger}$Department of Computer Science, University College London%
}

\begin{document}

\maketitle

\begin{abstract}
Unfolded proximal neural networks (PNNs) form a family of methods that combines deep learning and proximal optimization approaches. They consist in designing a neural network for a specific task by unrolling a proximal algorithm for a fixed number of iterations, where linearities can be learned from prior training procedure. PNNs have shown to be more robust than traditional deep learning approaches while reaching at least as good performances, in particular in computational imaging. However, training PNNs still depends on the efficiency of available training algorithms. In this work, we propose a lifted training formulation based on Bregman distances for unfolded PNNs. Leveraging the deterministic mini-batch block-coordinate forward-backward method, we design a bespoke computational strategy beyond traditional back-propagation methods for solving the resulting learning problem efficiently. We assess the behaviour of the proposed training approach for PNNs through numerical simulations on image denoising, considering a denoising PNN whose structure is based on dual proximal-gradient iterations. %We further propose a deterministic batch formulation of the learning strategy. 

\end{abstract}

\begin{keywords}
Image Gaussian denoising, unfolding, proximal neural networks, Bregman distance
\end{keywords}

\section{Introduction}

% Inverse problems are encountered broadly across a diverse range of imaging and sensing applications, spanning from areas in physics, astronomy, geophysics and medicine \cite{}. The intricacy of inverse problems often lies in formulating strategies to combat its ill-posedness \cite{}. This gives rise to interests in the development as a research field in  mathematics as well as applied fields \cite{}. 

In the last decade, deep learning approaches have become state-of-the-art in solving a great variety of tasks in data science \cite{Goodfellow-et-al-2016}. In particular, they have shown to be highly successful for solving image restoration problems. In this context, feedforward convolutional neural networks (CNNs) are commonly used, such as celebrated DnCNN \cite{zhang2017beyond}, Unet \cite{ronneberger2015u} and their variations. Such CNNs can be formulated as a composition of operators
\begin{equation} \label{def:nn}
    (\forall x \in \R^N)\quad
    G_{\thetab}(x) = F_{\theta_K} \circ T_{\theta_{K-1}} \cdots \circ T_{\theta_0}(x),
\end{equation}
where $\thetab = (\theta_k)_{1\le k \le K} \in \Theta$ are the underlying parameters of the network (e.g. convolution kernels, biases), $F_{\theta_K}\colon \R^{N_{K}} \to \R^{N}$ is an affine function parametrized by $\theta_K$
and, for every $k \in \{0, \ldots, K-1\}$,
\begin{equation} \label{def:lay}
    % (\forall k \in \{0, \ldots, K-1\})\;
    T_{\theta_k} \colon \R^{N_k} \to \R^{N_{k+1}} \colon u \mapsto D_k ( F_{\theta_k} (u) ),
\end{equation}
with $F_{\theta_k}\colon \R^{N_k} \to \R^{N_{k+1}}$ being an affine function parametrized by $\theta_k$, and $D_k \colon \R^{N_{k+1}} \to \R^{N_{k+1}}$ being an activation function.
More recently, inspired by the seminal work of Gregor and Lecun \cite{gregor2010learning} on LISTA, the optimization community started developing unfolded networks based on proximal optimization algorithms, namely proximal neural networks (PNNs). These network architectures basically unroll a fixed number of iterations of proximal algorithms (then renamed \textit{layers}), where different linear operators can be learned at each layer. Within the formulation \eqref{def:nn}-\eqref{def:lay}, for every $k\in \{0, \ldots, K\}$, $D_k$ corresponds to a proximity operator, and $F_{\theta_k}$ often corresponds to a gradient step. It has been shown that most of the activation functions used for feedforward networks correspond to proximity operators \cite{zhang2017convergence,hasannasab2020parseval,combettes2020deep}. %\cite{combettes2015?} 

On the one hand, PNNs have shown to outperform traditional proximal algorithms, and to reach similar performance as advanced CNNs for image restoration tasks \cite{adler2018learned,jiu2021deep}. On the other hand, since their structure is reminiscent of optimization theory, they also gain in interpretability and robustness \cite{monga2021algorithm}. 
Nevertheless, the robustness gained in the network architecture is still limited by the training procedure, as deep neural networks rely on highly parametrized nonlinear systems. 

Standard methods for learning parameters $\thetab$ (such as SGD and Adam \cite{kingma2014adam}) employ first order optimization methods where the (sub-)gradient information is evaluated using the back-propagation algorithm \cite{rumelhart1986learning}. Despite its popularity, this standard approach can suffer from potential drawbacks. 
First, the associated minimization problem is typically non-convex and existing algorithms offer no guarantee of optimality for the delivered output parameters. 
A second major issue is vanishing and exploding gradients, where gradients become extremely small or large, causing computation challenges, especially for deeper architectures \cite{he2016deep}. %While including skip-connections in network architectures may mitigate this issue, it still requires meticulous hyper-parameter tuning. 
Finally, the sequential nature of back-propagation makes it difficult to distribute computation for individual partial derivatives across multiple workers, hence constraining the overall training efficiency. 

These limitations have inspired fruitful research in seeking alternative training methods \cite{zach2019contrastive, zhang2017convergence,gu2020fenchel,taylor2016training,carreira2014distributed, xu2022alternative, frecon2022bregman, combettes2023variational, wang2023lifted}. One strategy consists in reformulating the original network parameter estimation problem as a penalized problem that can be solved by distributed optimization methods. In particular, in \cite{wang2023lifted} the lifted Bregman (LB) training strategy has been proposed, that splits the training process over the different layers and associates the proximal activation functions with tailored Bregman distance functions. By construction, the LB formulation leads to a relaxed problem, involving a collection of bi-convex penalty functions, paving the way to the use of advanced optimization algorithms for the parameter estimation. Due to its particular structure that leverages proximity operator theory, the LB approach appears to be particularly well suited for training PNNs. Hence, in this work, we explore a bespoke training procedure for PNNs in the simplified context of image denoising, leveraging the LB approach \cite{wang2023lifted}. 

%\textcolor{red}{In \cite{} a similar approach has been proposed, namely Bregman bla bla. It consists in [2 sentences!! say that (i) convexity properties? (ii) relies on proximity operator theory].}

The paper is structured as follows. 
We provide backgrounds on optimization and PNNs in Section~\ref{sec:background}. In Section~\ref{sec:learning-strategies} we review lifted training strategies, focusing on the LB formulation. 
We then introduce in Section~\ref{sec:proposed} a LB formulated training strategy for our PNN.
In Section~\ref{sec:numerics}, we give numerical experiment setups and show the efficiency of the proposed approach for image denoising. We conclude in Section~\ref{sec:conclusion}.

\section{Background}\label{sec:background}

\subsection{Notation}\label{sec:notation}

We begin this section by first introducing useful definitions. Let $\psi \colon \R^N \to ]-\infty, +\infty]$ be a proper, lower semi-continuous and convex function. 
% The Moreau sub-differential of $\psi$ at $v \in \R^N$ is defined as $\partial \psi(v) = \big\{ r \in \R^N \, | \, (\forall u \in \R^N)\; \langle u-v, r \rangle+ \psi(v) \le \psi(u) \big\}$.
The proximity operator of $\psi$ at $v \in \R^N$ is defined as $\prox{\psi}(v) = \argmind{u\in \R^N} \psi(u) + \frac12 \|u-v\|^2$.
The Fenchel-Legendre conjugate function of $\psi$ at $u\in \R^N$ is defined as $\psi^\ast(u) = \sup_{x \in \R^N} \left\{ \langle x, u \rangle - \psi(x) \right\}$.
% \begin{equation}
%     \partial \psi(v) = \Big\{ q \in \R^N \, | \, (\forall u \in \R^N)\; \langle u-v, q\rangle+ \psi(v) \le \psi(u) \Big\}.
% \end{equation}
% The generalized Bregman distance \cite{bregman1967relaxation, kiwiel1997proximal} of $\psi$, can be defined using the Fenchel-Legendre conjugate as %that measures the difference in function value between $\psi(u)$ for $u\in \R^N$ and its linearisation around point $v\in \R^N$, 
% \begin{equation}\label{def:bregman-distance}
% % (\forall (u,v) \in \R^N)
% % (\forall q(v) \in \partial \psi(v))\\
%     % D_\psi^{r(v)}(u, v) = \psi(u) - \psi(v) - \langle r(v), u - v \rangle ,
%     D_\psi^{r(v)}(u, v) = \psi(u) + \psi^\ast(r(v)) - \langle u, r(v) \rangle
% \end{equation}%
% where $r(v) \in \partial \psi(v)$ (see \cite[Theorem 23.5]{rockafellar1970convex}).
% Here, $q(v) \in \partial \psi(v)$ is a subgradient of $\psi$ at argument $v \in \mathbb{R}^n$. 
% The (generalised) Bregman distances measures the difference in function value between $\psi(u)$ and its linearisation around point $v$. For more examples and use cases we refer readers to \cite{censor1981iterative, burger2016bregman,benning2021bregman}. 
% For a proper, convex and semi-continuous function $F$, we denote its convex- or Fenchel-conjugate by $F^\ast$, defined as
% \begin{align}\label{[def:fenchel-conjugate]}
%     \psi^\ast(y) = \sup_{x} \left\{ \langle x, y \rangle - \psi(x) \right\} \, .
% \end{align}
% Then the Bregman distance can be expressed as $D_\psi^{r(v)}(u, v) = \psi(u) + \psi^\ast(r(v)) - \langle u, r(v) \rangle$ \cite[Theorem 23.5]{rockafellar1970convex}.
%
The Bregman penalty function \cite{wang2023lifted} associated with $\psi$ at $(u,v) \in (\R^N)^2$ is defined as
\begin{equation}  \label{def:bregman-loss}
    \bregmanloss(u, v) := q_{\psi}(u) + q_\psi^\ast\left( v \right) - \langle u, v \rangle  
\end{equation} 
where $q_{\psi}(u) = \frac12 \| u \|^2 + \psi(u)$
% \begin{equation}\label{eq:q-psi}
%     q_{\psi}(u) = \frac12 \| u \|^2 + \psi(u) ,
% \end{equation}
and $q_\psi^\ast$ refers to the Fenchel-Legendre conjugate of $q_\psi$ . 
% \begin{equation*}
%     \left( \frac12 \| \cdot \|^2 + \psi \right)^\ast (y) = \sup_{x} \; \langle y, x \rangle - \frac12 \|x\|^2 - \psi(x) \,.
% \end{equation*}
As a Bregman distance \cite{bregman1967relaxation, kiwiel1997proximal}, $ \bregmanloss(u, v)$ is non-negative and, for $v\in \R^N$ fixed, its global minimum value is zero and $\prox{\psi}(v) = \argmind{u \in \R^N} \bregmanloss(u, v)$. %is attained at $u = \prox{\psi}(v)$, i.e.
% \begin{equation}
%     \prox{\psi}(v) = \argmind{u \in \R^N} \bregmanloss(u, v).
% \end{equation}
It is differentiable with respect to its second argument, with
\begin{equation}\label{eq:grad-bregman-loss}
(\forall (u, v)\in (\R^N)^2)\quad
    \nabla_v \bregmanloss[\psi](u,v) = \prox{\psi}(v) - u  .
\end{equation}

%We now proceed to introduce the denoising problem. 

\subsection{Image denoising}\label{sec:image-denoising}

%In this work, we focus on a key type of inverse problem in the field of imaging, namely the image restoration problems. In contrary to image reconstruction problems, this type of inverse problems assumes that the observations are corrupted realization of the orignal images hence the data space and the image space coincide \cite{}. 
An image denoising problem seeks to find an estimate $\widehat{x}\in \R^N$ of an original clean image $\overline{x} \in \R^N$ from noisy observations $z \in \R^N$. We focus on the case of additive Gaussian noise, where the forward model is given by $z = \bar{x} + \epsilon$
% \begin{equation}
%     z = \bar{x} + \epsilon \,, \label{eq:denoising-problem}
% \end{equation}
where $\epsilon \in \R^N$ is a realization of a random zero-mean Gaussian variable with standard deviation $\sigma>0$. A common strategy consists in defining $\widehat x$ as % for finding estimates is via the penalised least-squares approach
\begin{align}
    \widehat{x} = \argmind{x \in \R^N} \frac{1}{2} \|x-z\|^2 +  g(Lx), \label{eq:denoise-objective} 
\end{align}%
where $g \colon \R^S \to \R^N$ is some regularization function and $L \colon \R^N \to \R^S$ is a linear operator \cite{rudin1992nonlinear, mallat1999wavelet, jacques2011panorama}. %, and $\lambda>0$ is a regularisation parameter. %Finally, the least-squares data-fidelity term and the regularisation function term are balanced by the associated non-negative regularisation parameter $\lambda>0$. 
% A popular choice for the regularisation function $g \circ L$ is the $\ell_1$ norm, chosen to promote sparsity in the transformed domain induced by $L$ \cite{rudin1992nonlinear, mallat1999wavelet, jacques2011panorama}. %add references here include [$\ell_1$+wavelet, TV, etc]. 

\subsection{PNNs for image denoising}\label{sec:PNNs}

% For decades, proximal algorithms have been the start-of-the-art approach for solving this type of variational problem due to the non-smoothness and composited nature of $g \circ L$ \cite{combettes2011proximal}. %Depends on specific problem structure, various types of splitting methods can be employed to computationally efficiently obtaining desirable solutions \cite{}.
% In particular, 
Problem~\eqref{eq:denoise-objective} can be solved efficiently by the dual forward backward algorithm (dual-FB) \cite{combettes2010dualization}, given by
\begin{equation}\label{algo:dfb-dual}
    \begin{array}{l}
        \text{for } k = 0, 1, \ldots \\
        \left\lfloor
        \begin{array}{l}
            \tau_k \in [\underline{\tau}, (2-\underline{\tau})\| L \|^{-2}] \\
            u_{k+1} = \prox{\tau_k g^*} \big( u_{k} - \tau_k L (L^* u_{k} - z)  \big),
            % u_{k+1} = \prox{\tau_k g^*} \Big( (\Idt - \tau_k L L^*) u_{k} + \tau_k L z  \Big),
        \end{array}
        \right.
    \end{array}
\end{equation}
where $\underline{\tau}>0$. Then, according to \cite{combettes2010dualization}, we have $\widehat{x} = \lim_{k\to +\infty}  z-L^* u_k$.
% \begin{equation}
%     \widehat{x} = \lim_{k\to +\infty}  z-L^* u_k.
% \end{equation}
% For appropriate choices of the step-sizes $(\tau_k)_{K\in \N}$, $(u_k)_{k\in \N}$ converges to a solution $\widehat u$ of the dual problem.  $L^*$ is the adjoint operator of $L$, and $  g^*$ denotes the Fenchel conjugate function of $ g$.Then a primal solution to~\eqref{eq:denoise-objective} is obtained as $\widehat x = z-L^* \widehat u$. 
%summarized as the following
%\begin{equation}\label{algo:dfb-const}
%    \begin{array}{l}
%        u_0 \in \R^S, \\
%        \text{for } k = 0, 1, \ldots \\
%        \left\lfloor
%        \begin{array}{l}
%            x_{k} = \text{Proj}_{C} (z - L^{*}u_{k}), \\
%            u_{k+1} = \text{prox}_{\tau_{k} g^{*}}(u_{k} + \tau_{k} L x_{k}),
%        \end{array}
%        \right.
%    \end{array}
%\end{equation}
%where, for every $k\in \N$, $\tau_k>0$ is a step-size chosen to ensure convergence of sequences $(x_k)_{k\in \N}$ and $(u_k)_{k\in \N}$, $L^*$ is the adjoint operator of $L$, and $  g^*$ denotes the Fenchel conjugate function of $ g$.

%le2023
Recently, a few works have proposed to unroll the dual-FB algorithm to design unfolded PNNs. In particular, in \cite{repetti2022dual, le2022faster, le2023pnn}, the authors unrolled the dual-FB algorithm for the denoising problem, and shown that despite having $10^3$ less parameters than a DRUnet \cite{zhang2021plug}, it reaches similar performances. We hence adopt a similar strategy and unroll algorithm~\eqref{algo:dfb-dual} over a fixed number of iterations $K \in \N^*$ to
construct a composited nonlinear mapping $G_{\thetab}$ of the form of~\eqref{def:nn}-\eqref{def:lay} with
% The resulting unfolded neural network $G_{\thetab}$ takes the form of a proximal neural network \cite{} with underlying parameters $\thetab = (\theta_k)_{0 \le k \le K} \in \Theta $, where the inner layers of the network are defined by
\begin{equation*}
\begin{cases}
    T_{\theta_0} \colon \R^N \to \R^S \colon x \mapsto L_{\theta_0} x \\
    F_{\theta_K} \colon \R^S \to \R^N \colon u \mapsto z - L_{\theta_K}^* u,
\end{cases}
\end{equation*}
% \begin{align}
%     &T_{\theta_0} \colon \R^N \to \R^S \colon x \mapsto L_{\theta_0} x \\
%     &F_{\theta_K} \colon \R^S \to \R^N \colon u \mapsto z - L_{\theta_K}^* u.
% \end{align}
and for every $k \in \{1, \ldots, K-1\}$,
\begin{equation} \label{def:PNN-layer-gen}
(\forall u \in \R^S)\quad
    T_{\theta_k} = \widetilde{T}_{\theta_k}^{(\mathcal{D})} \Big( \widetilde{T}_{\theta_{k}}^{(\mathcal{P})}(u) \Big)
\end{equation}
with
% \begin{equation}
%     \begin{array}{l@{}l@{}c@{}l}
%         T_{\theta_{k}} \colon
%         &   \R^S &\to & \,\R^S \\
%         &   u &\mapsto & \,\prox{\tau_k g^*} \Big( (\Idt - \tau_k L_{\theta_k} L_{\theta_k}^*) u + \tau_k L_{\theta_k} z  \Big).
%     \end{array}
% \end{equation}
\begin{align}
    & \widetilde{T}_{\theta_{k}}^{(\mathcal{P})} \! \colon \R^S \!\! \to  \R^N\!\! \times\! \R^S \colon u \mapsto  \big( L_{\theta_k}^* u +  z , \, u \big) , \label{def:PNN-layer-P}\\
    & \widetilde{T}_{\theta_k}^{(\mathcal{D})}\! \colon \R^N\!\!\! \times\! \R^S \!\!\! \to \! \R^N \! \colon \!(x,u) \! \mapsto  \prox{\tau_k g^*} \! ( u - \tau_k L_{\theta_k} x ) \label{def:PNN-layer-D}
\end{align}
and $\tau_k \in [\underline{\tau}, (2-\underline{\tau}) \| L_{\theta_k} \|^{-2}]$.
Each intermediate layer is a composition of a proximal activation function over an affine transformation where the underlying parameters of $(T_{\theta_k})_{0 \le k \le K-1}$ and $F_{\theta_K}$ are $\thetab = (\theta_k)_{0 \le k \le K} $, corresponding to a parametrization of the linear operators $(L_{\theta_k})_{0 \le k \le K}$. We define by $\Theta$ the feasible space for $\thetab$. In this context, the unfolded network allows the use of different linear operator $L_{\theta_k}$ across layers, hence providing higher adaptivity to the task of interest by enlarging the feasibility space. 
%{\color{red}TO ADD: quadratic formulation -- not affine -- for combined primal-dual layers}

\section{Learning strategies}\label{sec:learning-strategies}

% We begin this section by first defining the backbone neural network architecture in this work, namely the dual-FB net. We next proceed to the discussion of the lifted Bregman training strategy for learning neural network parameters. 

\subsection{Lifted training strategies}

In the context of supervised learning for image denoising, we assume that we have access to a training dataset composed of couples of groundtruth/noisy images $(\overline{x}^{(s)}, z^{(s)})_{ s \in \mathbb S} $. Then, standard supervised learning approaches often aim to
\begin{equation}\label{pb:erm}
    \minimize{\thetab \in \Theta} \dfrac{1}{\#(\eS)} \sum_{s \in \eS} \ell \Big( G_{\thetab}(z^{(s)}), \overline{x}^{(s)} \Big) \,,
\end{equation}
where $\ell \colon \R^N \times \R^N \to \R$ is a data error function (e.g., $\ell^2$ loss) measuring discrepancy between the output of the network $G_{\thetab}(z)$ and the groundtruth $\overline{x}$. 
%The training dataset of $\#(\eS)$ pairs of clean and noisy samples is denoted by $\eS = \{\overline{x}^{(s)},z^{(s)}\}_{s \in \eS}$. 
Standard computational approaches for solving \eqref{pb:erm} are (sub-)gradient based algorithms (e.g., (sub-)gradient descent or its stochastic variants), where the evaluation of the gradient with respect to $\thetab$ is computed using back-propagation \cite{rumelhart1986learning}.

In the context of feedforward networks of the form of~\eqref{def:nn}-\eqref{def:lay}, problem~\eqref{pb:erm} can be split over layers by introducing auxiliary dual variables such as, for every $s \in \mathbb S$, $u_0^{(s)} = T_{\theta_0}(z^{(s)})$, $\overline{x}^{(s)} = F_{\theta_K}(u_{K-1}^{(s)})$, 
% \begin{equation}
%     u_0^{(s)} = T_{\theta_0}(z^{(s)}), 
%     \quad
%     \overline{x}^{(s)} = T_{\theta_K}(u_{K-1})
% \end{equation}
and, for every $k \in \{1, \ldots, K-1\}$, $u^{(s)}_{k} = T_{\theta_k}(u^{(s)}_{k-1})$.
% \begin{equation}
%     u^{(s)}_{k} = T_{\theta_k}(u^{(s)}_{k-1}) .
% \end{equation}
%
%
% In the particular case of the PNN presented in Section~\ref{sec:PNNs}, %it can be noticed that %each pair of consecutive iterates in the dual-FB algorithm are associated by the following relation:
% we can introduce
% \begin{equation}
%     u_0^{(s)} = T_{\theta_0}(z^{(s)}), 
%     \quad
%     \overline{x}^{(s)} = T_{\theta_K}(u_{K-1})
% \end{equation}
% \begin{align*}
%     &   u_0^{(s)} = T_{\theta_0}(z^{(s)})  \\
%     (\forall k \in \{1, \ldots, K-1\})\quad 
%     &   u^{(s)}_{k} = T_{\theta_k}(u^{(s)}_{k-1}) .
% \end{align*}
% To ease the notation, we introduce $\ub = (u_k)_{0 \le k \le K-1} \in \mathcal{U} $ to denote the collection of all auxiliary dual variables. 
Hence, Problem~\eqref{pb:erm} can equivalently be written as a constrained minimization problem:
\begin{multline} \label{pb:min-loss-constrained}
    \minimize{\thetab , (\ub^{(s)})_{s \in \mathbb S}} \dfrac{1}{\#(\eS)} \sum_{s \in \eS} \ell \Big( F_{\theta_K}(u_{K-1}^{(s)}), \overline{x}^{(s)} \Big)  \\
    \text{s.t. }
    \begin{cases}
    u_0^{(s)} = T_{\theta_0}(z^{(s)}) \\
    (\forall k \in \{1, \ldots, K-1\})\;
    u^{(s)}_{k} = T_{\theta_k}(u^{(s)}_{k-1}) ,
    \end{cases}
\end{multline}%
where, for every $s \in \mathbb S$, $\ub^{(s)} = (u_k^{(s)})_{0 \le k \le K-1} $ denotes the auxiliary dual variables for a given pair $(\overline{x}^{(s)}, z^{(s)})$.
In this equivalent formulation, the popular back-propagation algorithm can be deduced from a Lagrangian formulation \cite{lecun1988theoretical}. In this context, alternative pathways to standard stochastic back-propagation algorithms can then be considered. Nevertheless, \eqref{pb:min-loss-constrained} necessitates to handle a collection of non-linear constraints, which can be inefficient in practice.

An approach to overcome this issue is to adopt a penalty approach, where a relaxed formulation of~\eqref{pb:min-loss-constrained} is considered by replacing the constraints by penalty terms \cite{carreira2014distributed,taylor2016training}. In the next section, we review one such approach that appears to be particularly suited for training PNNs.
% discuss a more dedicated approach, the LB training framework, \textcolor{red}{[weak:] whose properties enjoy many computational convenience + prox theory + adapted to PNNs bla bla.} 

\subsection{Lifted Bregman approach}
\label{sec:LB-back}

Following \cite{wang2023lifted}, a relaxed version of problem~\eqref{pb:min-loss-constrained} consists of penalizing the constraints in~\eqref{pb:min-loss-constrained} with penalty functions instead of strictly enforcing them. Here, the penalty functions are chosen as a compositions of the Bregman functions defined in \eqref{def:bregman-loss} and the affine-linear transformations $F_{\theta_k}$. More precisely, we consider the particular case where activation functions $(D_k)_{0 \le k \le K-1}$ in \eqref{def:lay} are proximity operators of some functions $(\psi_k)_{0 \le k \le K-1}$, and $N_0 = N$ and $N_1 = \ldots = N_{K} = S$ (i.e., all inner-layers have input/output of the same dimension).
This leads to the relaxed formulation %In this context, the training objective is lifted by the Bregman penalty functions, defined as follows. 
\begin{multline}\label{pb:min-loss-Breg-gen}
    \minimize{\thetab , (\ub^{(s)})_{s \in \mathbb S}} \dfrac{1}{\#(\eS)} \bigg(
    \sum_{s \in \eS} \ell \big( F_{\theta_K}(u_{K-1}^{(s)}), \overline{x}^{(s)} \big)   \\[-0.2cm]
   \!+\! \bregmanloss[\psi_0](u^{(s)}_0 \!, F_{\theta_0}(z^{(s)})) \!+ \! \sum_{k=1}^{K-1} \! \bregmanloss[\psi_k](u^{(s)}_k, F_{\theta_k}(u_{k-1}^{(s)})) \bigg).
\end{multline}%
%\noindent where $\bregmanloss$ is the Bregman penalty defined in~\eqref{def:bregman-loss}. 
The LB approach can be viewed as a generalization of some classical lifted training methods \cite{ carreira2014distributed, taylor2016training}. However, unlike quadratic penalty approaches that still need to differentiate non-smooth activation functions when first-order optimization methods are applied, the Bregman penalty function $\bregmanloss$ is continuously differentiable with respect to its second argument (see \eqref{eq:grad-bregman-loss}), and hence with respect to the network parameters $\thetab$.
%
% \textcolor{red}{[Use notation in \eqref{def:lay} to explain how the constraint in \eqref{pb:min-loss-constrained} are relaxed (give equations !)]}
% For a proper, lower-semicontinuous and convex function $\psi$, for all $(x,y) \in \R^N \times \R^N $, the associated Bregman penalty function is a Bregman distance function defined as 
% \begin{equation}  
%     \bregmanloss(x, y) := q_{\psi}(x) + q_\psi^\ast\left( y \right) - \langle x, y \rangle \, ,\label{eq:bregman-loss}
% \end{equation} 
% where
% \begin{equation}\label{eq:q-psi}
%     q_{\psi}(x) = \frac12 \| x \|^2 + \psi(x) \,.
% \end{equation}
% The notation $q_\psi^\ast$ refers to the Fenchel conjugate of $q_\psi$, i.e., 
% \begin{equation*}
%     \left( \frac12 \| \cdot \|^2 + \psi \right)^\ast (y) = \sup_{x} \; \langle y, x \rangle - \frac12 \|x\|^2 - \psi(x) \,.
% \end{equation*}
% The Bregman penalty function are non-negative and attains its global minimum value zero at $x = \prox{\psi}(y)$. 
% This property links a chosen proximal activation function to its associated Bregman penalty function. In other words, for every $y\in \R^N$, 
%     \begin{equation*}
%         \prox{\psi}(y) = \argmind{x \in \R^N} \bregmanloss(x, y).
%     \end{equation*}
%
Further, the following properties can be deduced from \cite[Theorem 10]{wang2023lifted} (see \cite{wang2023lifted, wang2023inversion} for more properties).

\vspace{-0.1cm}
\begin{proposition}\label{prop}
\begin{enumerate}
    \item\label{prop:i} Problem~\eqref{pb:min-loss-constrained} and problem~\eqref{pb:min-loss-Breg-gen} share the same set of solutions.

    \vspace{-0.2cm}
    \item\label{prop:ii} Problem~\eqref{pb:min-loss-Breg-gen} is convex with respect to $\thetab$ when $(\ub^{(s)})_{s \in \mathbb S}$ is fixed.

    \vspace{-0.2cm}
    \item\label{prop:iii} 
    For every $k\in \{0, \ldots, K-1\}$, problem~\eqref{pb:min-loss-Breg-gen} is convex with respect to $(\ub_k^{(s)})_{s \in \mathbb S}$ when $\thetab$ and the other layer's auxiliary variables are fixed.
\end{enumerate}    
\end{proposition}

\vspace{-0.1cm}

The LB training framework fits naturally to the structure of PNNs whose activation functions are proximity operators. In this context, computing the gradient of $\bregmanloss$ with respect to $\thetab$ would boil down to using the activation function itself. In the next section, we formulate the LB training problem for the PNN described in Section~\ref{sec:PNNs}. % and discuss computational strategies for efficiently solving the associated minimisation problem. 

%mention the biconvex formulation in the later computational properties section

%mention the biconvex formulation in the later computational properties section

\section{Proposed LB training for PNNs}
\label{sec:proposed}

In this section we design a LB training strategy for the denoising PNN described in Section~\ref{sec:PNNs}. % apply the lifted Bregman training strategy and formulate the lifted training objective in particular for the dual-FB network. 

\subsection{Proposed lifted Bregman formulation}
\label{ssec:prop:LB}

One can note that the inner layers of the PNN network described in Section~\ref{sec:PNNs} are themselves compositions of two feed-forward layers of the form of~\eqref{def:lay}, where the first sub-layers $\widetilde{T}^{(\mathcal{P})}$ (see \eqref{def:PNN-layer-P}) provide outputs in $\R^N \times \R^S$ and the second sub-layers $\widetilde{T}^{(\mathcal{D})}$ (see \eqref{def:PNN-layer-D}) take inputs from $\R^N \times \R^S$. Due to this change of space in the sub-layers, the LB formulation described in Section~\ref{sec:LB-back} cannot be directly applied, and instead we will need to apply it to the global inner-layers $(T_{\theta_k})_{1 \le k \le K-1}$ defined in \eqref{def:PNN-layer-gen}.

We thus introduce the following LB formulation for our PNN. For the training, we aim to
% Following the definition of the Bregman penalty functions, we relax the constraints in Eq.~\eqref{pb:min-loss-constrained} with a collection of penalty terms added to the objective function. In this context, we now aim to solve the following unconstrained minimisation problem
\begin{equation}    \label{pb:min-bregman}
    \minimize{\thetab, (\ub^{(s)})_{ s \in \eS}} \dfrac{1}{\#(\eS)} \sum_{s \in \eS} E \Big( \thetab, \ub^{(s)} \Big| z^{(s)}, \overline{x}^{(s)} \Big)
\end{equation}
where, for any $(\ub, \thetab) \in (\R^S)^{K-1} \times \Theta$ and $(\overline{x},z)$, we have
\begin{multline} \label{eq:def-E}
    E \big( \thetab, \ub \big| z, \overline{x} \big) 
       = \ell  \big( \widetilde{F}_{\theta_K} ( u_{K-1} , z ), \overline{x} \big)
    + \bregmanloss[\tau_1 g^*] \big( u_1, \widetilde{F}_{\theta_0, \theta_1} (  z ) \big)  \\
    + \sum_{k=2}^{K-1} \bregmanloss[\tau_k g^*] \big( u_{k}, \widetilde{F}_{\theta_k} ( u_{k-1}, z ) \big) 
\end{multline}
with
\begin{equation}\label{eq:def-E-lay}
\begin{cases}
    \displaystyle \widetilde{F}_{\theta_0, \theta_1} (  z ) 
        = (\Idt - \tau_1 L_{\theta_1} L_{\theta_1}^*) L_{\theta_0} z + \tau_1 L_{\theta_1} z , \\
        % = (\Idt - \tau_1 L_{\theta_1} L_{\theta_1}^*) L_{\theta_0} z + \tau_1 L_{\theta_1} z , \\
    \displaystyle \widetilde{F}_{\theta_k} ( u_{k-1}, z ) 
        % = (\Idt - \tau_k L_{\theta_k} L_{\theta_k}^*) u_{k-1} + \tau_k L_{\theta_k} z , \\
        = u_{k-1} - \tau_k L_{\theta_k} (L_{\theta_k}^* u_{k-1} - z) , \\
    \widetilde{F}_{\theta_K} ( u_{K-1},  z ) 
        = z- L_{\theta_K}^* u_{K-1}.
\end{cases}
\end{equation}
In this context, the loss~\eqref{eq:def-E} is only convex with respect to the auxiliary variables $(\ub^{(s)})_{s \in \mathbb S}$, but not with respect to the parameters $\thetab$, due to the quadratic terms appearing in \eqref{eq:def-E-lay}. Nevertheless, Proposition~\ref{prop}\ref{prop:ii} still holds in this case.

% \review{Due to the form of $F_k$, it is however not fully convex with respect to $\thetab$ as we only have linearity with respect to the input layer parameter $\theta_0$ and output layer parameter $\theta_K$. To obtain a fully bi-convex formulation on $(\thetab, \ub)$ we could further split all the layers between primal and dual space by introducing a collection of auxiliary primal variables.}

%\textcolor{red}{Due to the form of $F_1$, it is however not convex with respect to $\thetab$. To obtain a fully bi-convex formulation on $(\thetab, \ub)$ we could further split $F_1$ introducing another auxiliary variable $u_0 = L_0 z$. }

\subsection{Proposed computational strategy}

Solving~\eqref{pb:min-bregman}-\eqref{eq:def-E} requires to estimate the network parameters and the auxiliary dual variables. It can be noticed that \eqref{eq:def-E} can be rewritten as the sum of a Lipschitz-differentiable function and a proximable function, i.e., for every $\ub \in (\R^S)^{K-1}$, $\thetab \in \Theta$ and $(\overline{x}, z)\in (\R^N)^2$, 
%$E\big( \thetab, \ub \big| z, \overline{x} \big) =  h(\ub) + f(\ub, \thetab | z, \overline{x})$
\begin{equation}\label{eq:E-split}
% (\forall \ub \in (\R^S)^{K-1})\quad
    E\big( \thetab, \ub \big| z, \overline{x} \big) 
    =  h(\ub) + f(\ub, \thetab | z, \overline{x})
\end{equation}
where $h$ is a proximable lower-semicontinuous, proper, convex function, and $f$ is the Lipschitz-differentiable function. 
In particular, we define $h(\ub) = \sum_{k=1}^{K-1} \tau_k g^* \big( u_{k} \big)$, while $f$ contains all the other terms.
% \begin{equation*}
%     h(\ub) = \sum_{k=1}^{K-1} \tau_k g^* \big( u_{k} \big),
% \end{equation*}
% and
% \begin{multline*}
%     f(\ub, \thetab | z, \overline{x}) 
%     =  \ell \big( \widetilde{F}_{\theta_K}(u_{K-1}, z), \overline{x} \big) 
%         + f_{\theta_0, \theta_1}(u_1,  z )  \\
%     + \sum_{k=2}^{K-1} f_{\theta_k}(u_{k-1}, u_{k},  z) ,
% \end{multline*}
% where, 
% \begin{multline*}
%     f_{\theta_0, \theta_1}(u_1,  z ) =  \frac{1}{2}\|u_1\|^2 + q_{\tau_1 g^*}^* \big( \widetilde{F}_{\theta_0, \theta_1}\big( z \big) \big) \\ 
%     - \left\langle u_{1}, \widetilde{F}_{\theta_0,\theta_1} \big( z \big) \right\rangle \,. 
% \end{multline*}
% and, for every $k\in \{1, \ldots, K-1\}$,
% \begin{multline*}
%     f_{\theta_k}(u_{k-1}, u_{k},  z) 
%     = \frac12 \|u_k\|^2 + q_{\tau_k g^*}^* \big( \widetilde{F}_{\theta_k}\big(u_{k-1}, z \big) \big) \\ 
%     - \left\langle u_{k}, \widetilde{F}_{\theta_k} \big( u_{k-1}, z \big) \right\rangle \, ,
% \end{multline*}

Then, we can rely on a block-coordinate forward-backward (FB) strategy to solve~\eqref{pb:min-bregman}, alternating between the optimization of $\thetab$ and the optimization of the auxiliary variables, whose iterations can be summarized as
%\textcolor{red}{[What? do we alternate over $\thetab$ and $\ub$, or over all the layers as well?]} 
% For a fixed pair of data $(z,\overline{x})$, the minimisation problem can be solved using the proximal gradient descent algorithm \cite{}, which results in the following iterative updates alternating between each $\theta_k$ and $u_k$
\begin{equation}\label{algo:BC-FB}
    \begin{array}{l}
        \text{for } n = 0, 1, \ldots \\
        \left\lfloor
        \begin{array}{l}
            \thetab^{(n + 1)} = \thetab^{(n)} - \beta_{n} \nabla_{\thetab} f(\ub^{(n)}, \thetab^{(n)} | z) \, ,\\
            \ub^{(n+1)} = \prox{\gamma_{n} \boldsymbol{\tau} g^*} \left( \ub^{(n)} - \gamma_{n} \nabla_{\ub} f(\ub^{(n)}, \thetab^{(n+1)} | z) \right)
        \end{array}
        \right.
    \end{array}
\end{equation}
% \begin{align}
%     \thetab^{(n + 1)} &= \thetab^n - \beta_{n} \nabla_{\thetab} f(\ub^n, \thetab^n | z) \, ,\\
%     \ub^{(n+1)} &= \text{prox}_{\gamma_{n} \boldsymbol{\tau} g^*} \left( \ub^{n} - \gamma_{n} \nabla_{\ub} f(\ub^n, \thetab^{(n+1)} | z) \right) \,,
% \end{align}
where $\boldsymbol{\tau} = \operatorname{Diag}((\tau_k)_{1 \le k \le K-1})$ and, for every $n\in \N$, $(\beta_{n} ,\gamma_{n})\in ]0,+\infty[^2$ are step-sizes that must be chosen to ensure convergence of the iterates generated by \eqref{algo:BC-FB} (see, e.g., \cite{chouzenoux2016block}). 
The partial gradients of $f$ in \eqref{algo:BC-FB} are evaluated using auto-differentiation combined with the differentiation property~\eqref{eq:grad-bregman-loss} of the Bregman function.
Finally, to use the proposed approach in a training setting, we consider a mini-batch approach in our simulations, where the partial gradients are approximated on sub-parts of the training dataset.
%We denote by $\boldsymbol{\tau}$ the collection of $\tau_k$ for $k \in \{1, K-1\}$.

\smallskip\noindent\textbf{Step-size backtracking strategy. \quad}
The choice of step-sizes $(\beta_n)_{n\in \N}$ and $(\gamma_n)_{n\in \N}$ in \eqref{algo:BC-FB} relies on the computation of the Lipschitz constants associated with the gradients of function $f$. In the considered training context, such a computation can be very expensive, and instead we propose a backtracking strategy %to ensure monotonicity of algorithm~\eqref{algo:BC-FB}.
% To avoid the expensive computations of the Lipschitz constants at each iteration of the training strategy, we propose a backtracking strategy to estimate $(\beta_n)_{n\in \N}$ and $(\gamma_n)_{n\in \N}$, 
that relies on monotone properties of the block-coordinate FB algorithm (that holds even for non-convex objective functions). 
Indeed, according to~\cite{chouzenoux2016block}, the step-sizes should be chosen to satisfy, for every $n\in \N$,
\begin{multline}\label{eq:maj-theta}
    f(\ub^{(n)}, \thetab^{(n+1)} | z)
    % \le \overline{f}_{\beta_n}(\ub^{(n+1)}, \ub^{(n)}, \thetab^{(n)} | z),
    \le f(\ub^{(n)}, \thetab^{(n)} | z) \\
        + \langle \thetab^{(n+1)} - \thetab^{(n)}, \nabla_{\thetab} f(\ub^{(n)}, \thetab^{(n)} | z)\rangle \\
        + \frac{1}{2 \beta_n} \| \thetab^{(n+1)} - \thetab^{(n)} \|^2
\end{multline}
and
\begin{multline}\label{eq:maj-u}
    f(\ub^{(n+1)}, \thetab^{(n+1)} | z)
    % \le \overline{f}_{\beta_n}(\ub^{(n+1)}, \ub^{(n)}, \thetab^{(n)} | z),
    \le f(\ub^{(n)}, \thetab^{(n+1)} | z) \\
        + \langle \ub^{(n+1)} - \ub^{(n)}, \nabla_{\ub} f(\ub^{(n)}, \thetab^{(n+1)} | z)\rangle \\
        + \frac{1}{2 \gamma_n} \| \ub^{(n+1)} - \ub^{(n)} \|^2.
\end{multline}
Since the partial gradients in \eqref{eq:maj-theta} and \eqref{eq:maj-u} need to be evaluated in the iterations~\eqref{algo:BC-FB}, verifying these inequalities appears to have a very low computational cost, that only necessitates extra forward passes to evaluate function $f$. We hence propose to use them to design our backtracking strategy. 
In particular, since the objective is non-convex, we adopt a \textit{two-way} backtracking approach (i.e., we automatically either increase or decrease the step-sizes at each iteration) in order to choose the largest step-sizes satisfying the majorant properties.

\section{Simulations and Results}\label{sec:numerics}

In this section, we aim to demonstrate the performance of the proposed strategy on an image denoising example.  
All numerical results are computed using Pytorch 2.2.2 on NVIDIA GeForce RTX 2080 Ti.

% Dual-FB networks are trained on subsets of the ImageNet dataset \cite{russakovsky2015imagenet} and evaluated on BSDS dataset \cite{martin2001database}. All images are pre-processed and resized to the size of $128 \times 128$. We consider two numerical experiments to examine the proposed strategy's stability and convergence performance in solving the denoising problem. Noisy observations are generated by $z^{(s)} = \bar{x}^{(s)} + \epsilon^{(s)}$ with $\epsilon^{(s)} \in \R^N$ being a realisation of a zero-mean Gaussian variable with standard deviation $\sigma = 0.1$. 

\subsection{Experiment setup}

\noindent\textbf{Datasets.\quad}
For the training groundtruth set, we consider a subset of $1,000$ images randomly selected from the Tiny ImageNet dataset \cite{russakovsky2015imagenet}, and we use patches of size $32 \times 32$. The associated noisy images are generated as, for every $s \in \mathbb S$, $z^{(s)} = \bar{x}^{(s)} + \epsilon^{(s)}$, with $\epsilon^{(s)} \in \R^N$ being a realization of a zero-mean Gaussian variable with standard deviation $\sigma = 0.1$. 
We then evaluate the resulting networks on $9$ images selected from the BSDS dataset \cite{martin2001database}.

\smallskip\noindent\textbf{Network architectures.\quad}
We evaluate our training procedure for the PNN defined in Section~\ref{sec:PNNs} for different network depth $K \in \{5, 10, 15\}$. For the three cases, we fix the number of features to $16$. We set $\tau_k$ at $1.8/\|L_k\|^2$, where the spectral norm is computed via power iteration following the method proposed in~\cite{repetti2022dual}. 
For all our experiments, the network weights are initialized following the strategy proposed in \cite{glorot2010understanding}.

\smallskip\noindent\textbf{Evaluation.\quad}
We consider two experiment settings to evaluate the performances of the proposed LB-FB training strategy. For both we use the $\ell^2$ loss and run algorithms over $50$ epochs. \\
\textit{Mini-batch stability:}
We fix the number of layers to $K=5$ and investigate the stability of the proposed alternating FB algorithm when considering a mini-batch approach. We vary the batch size at $\{1, 2, 5, 10, 20, 50, 100\}$ and train the networks on a small collection of $100$ images. \\
\textit{Comparison with SGD:}
We fix the mini-batch size to saturate memory, vary $K\in \{5, 10, 15\}$, and compare the proposed LB-FB training approach to a standard stochastic gradient descent approach (SGD) implemented in Pytorch. We run both algorithms on the $1,000$ images. For LB-FB, learning rates $(\beta_n, \gamma_n)_{n \in \N}$ are determined via the backtracking strategy. For SGD, the choice of learning rate is made via a grid search over values in the range $[1\times10^{-6},  5\times10^{-4}]$, where it is chosen empirically such that the training loss value is the lowest.

\begin{figure}[!t]
    \centering
    \includegraphics[width=0.49\columnwidth]{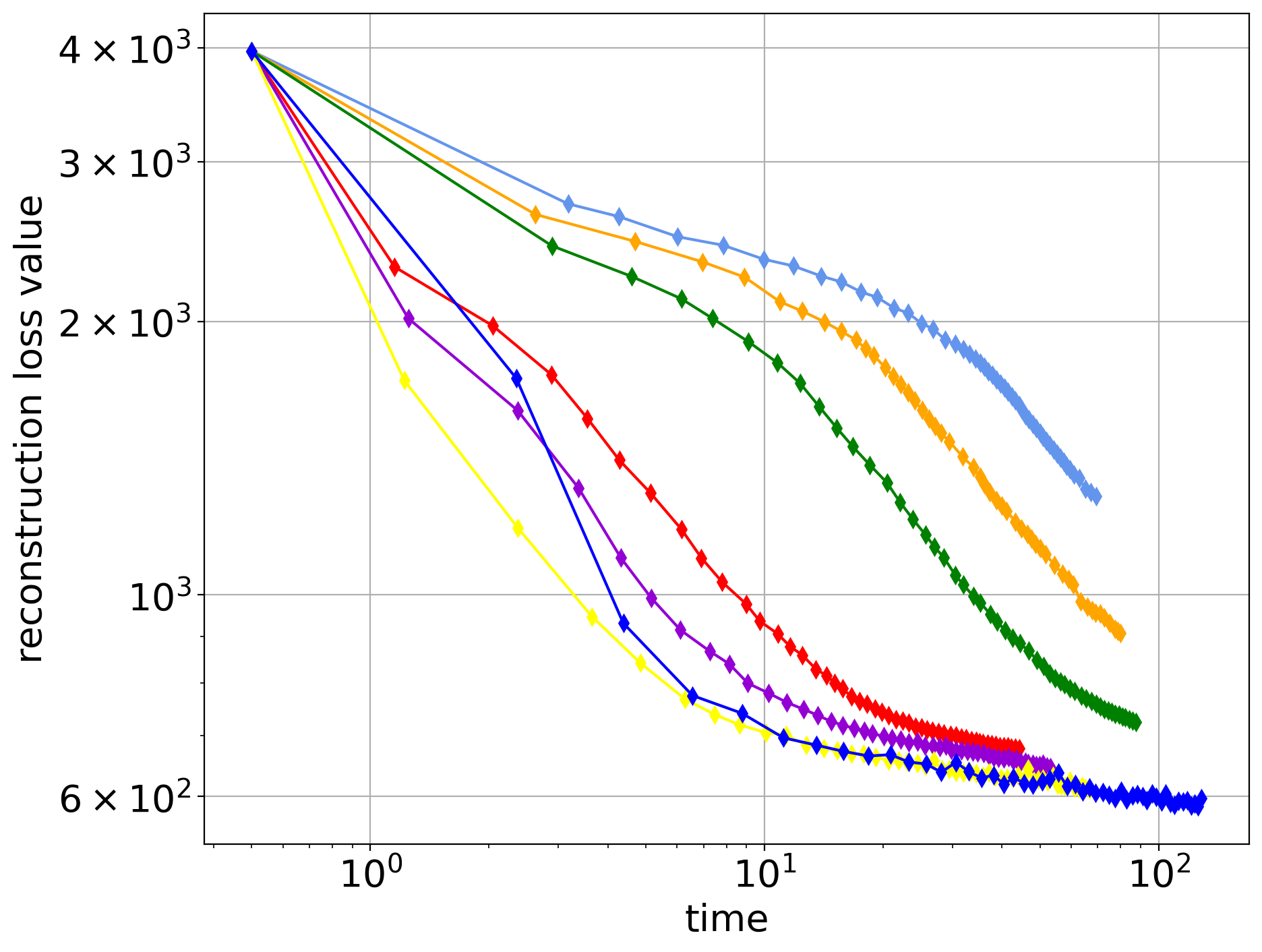}
    \includegraphics[width=0.49\columnwidth]{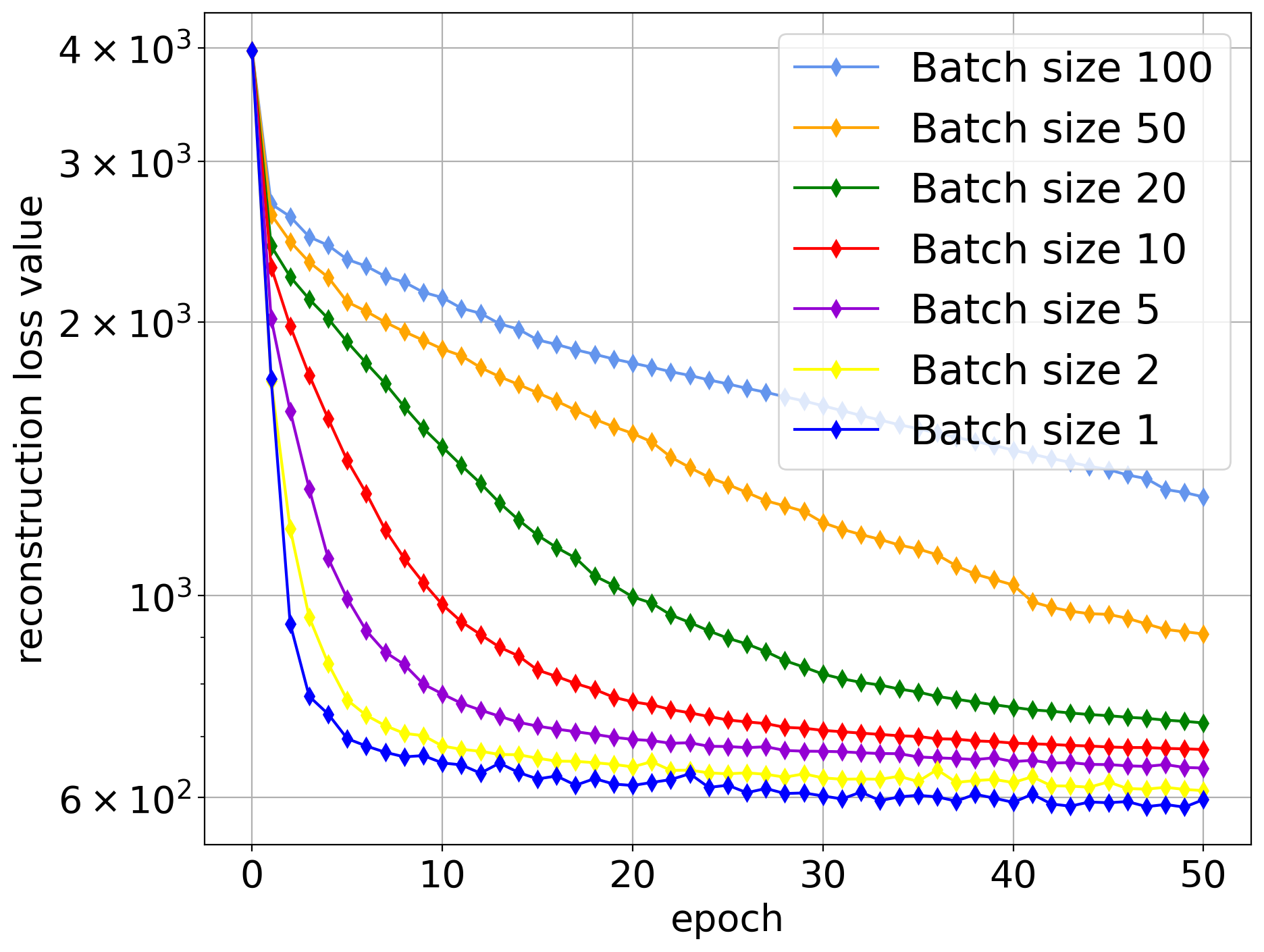}

    \vspace{-0.3cm}
    \caption{
    \small \textit{Mini-batch stability:} $\ell^2$ loss values with respect to computation time (left) and epochs (right), for different batch sizes. 
    % Examine proposed algorithm stability through varying mini-batch sizes. \textbf{Left:} loss value decay over run time. \textbf{Right: } loss value decay over epochs. 
    }
    \label{fig:numerics_batches}

% \vspace{-0.2cm}
\end{figure}

\begin{figure}[!t]
    \centering
    \includegraphics[width=0.49\columnwidth]{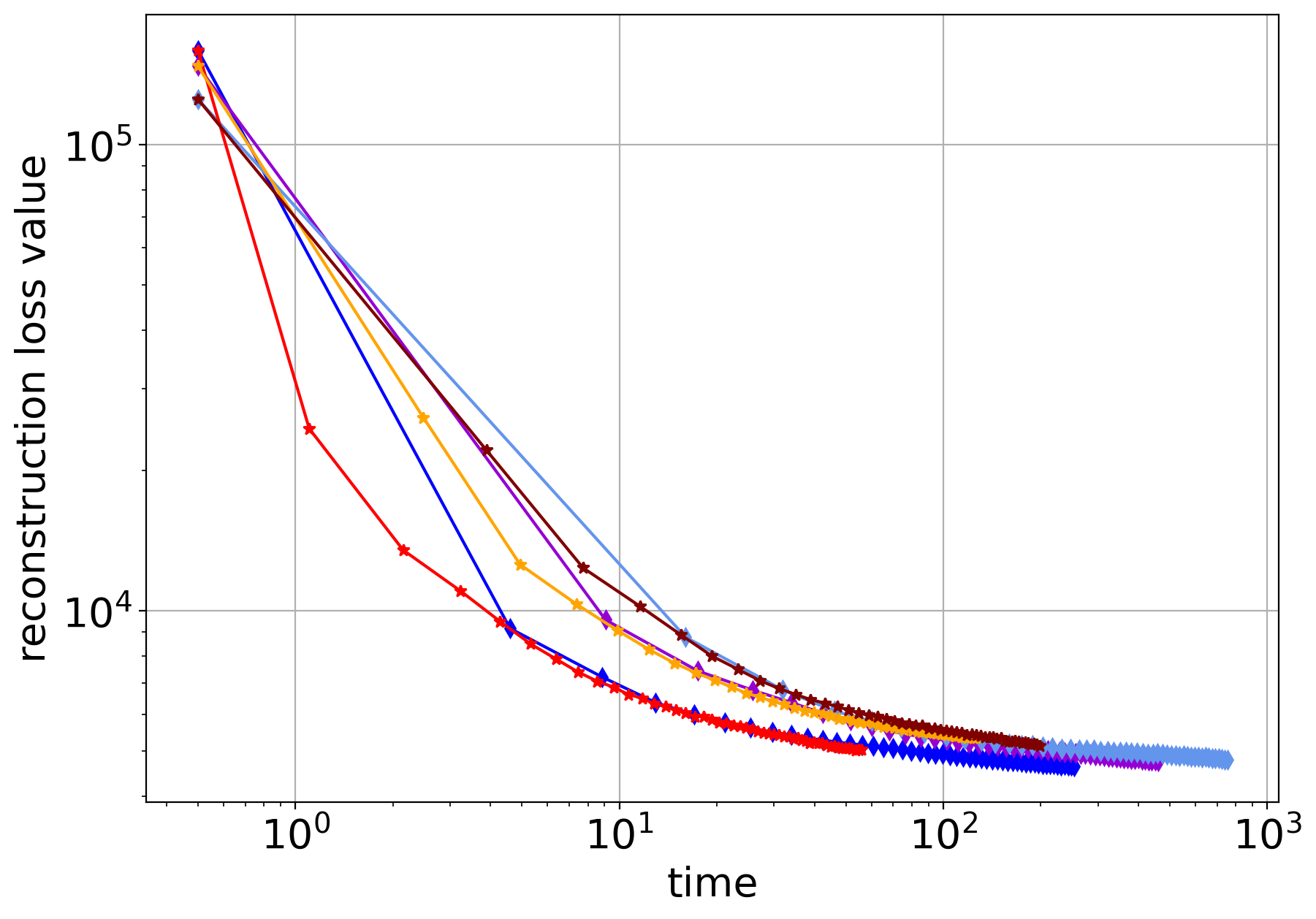}
    \includegraphics[width=0.49\columnwidth]{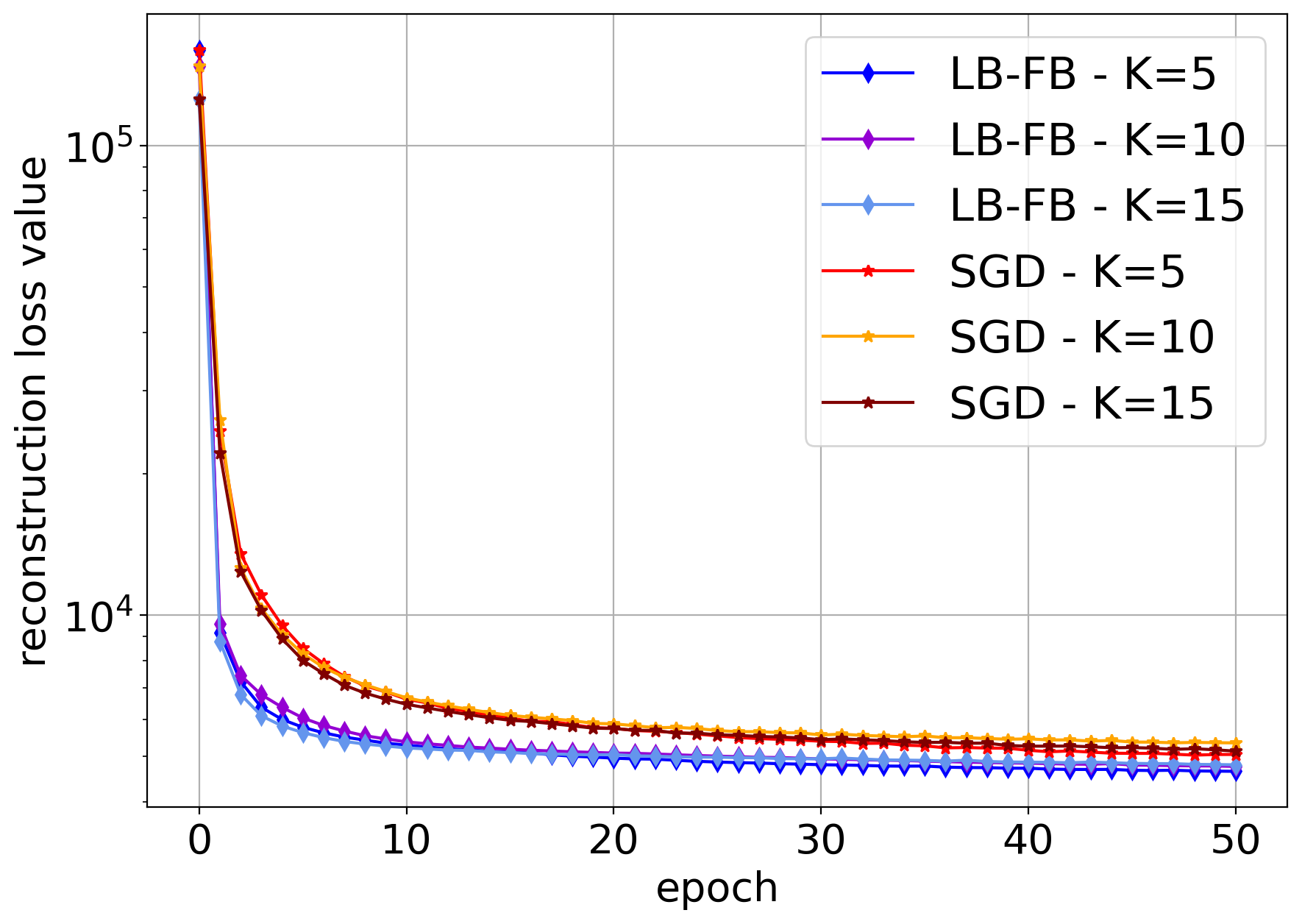}
    
    \vspace{-0.3cm}
    \caption{\small
    \textit{Comparison with SGD:} $\ell^2$ loss values comparing with SGD, with respect to computation time (left), and epochs (right).}
    \label{fig:numerics_convergence}

% \vspace{-0.2cm}
\end{figure}

\begin{figure}[!t]
\centering
\setlength\tabcolsep{0.0cm}
\begin{tabular}{ccc}
\includegraphics[scale=0.62, trim={0 0.5cm 0.25cm 0.25cm},clip]{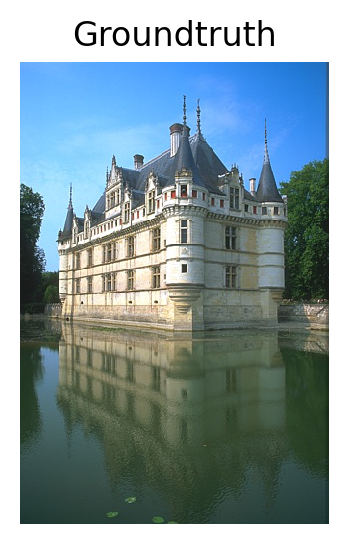}
&\includegraphics[scale=0.62, trim={0 0.5cm 0.25cm 0.25cm},clip]{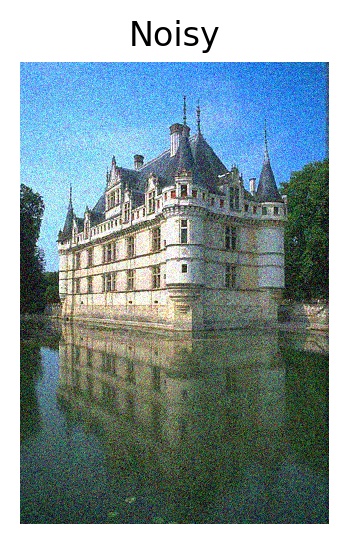}
&\includegraphics[scale=0.62, trim={0 0.5cm 0.25cm 0.25cm},clip]{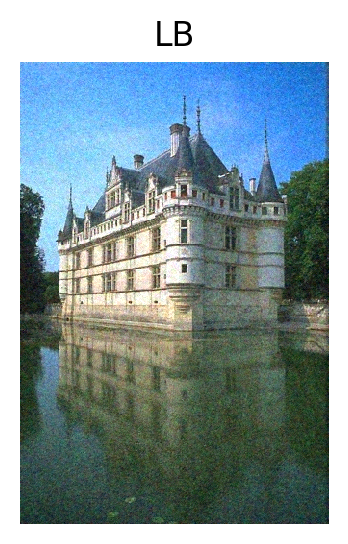} \\
\includegraphics[scale=0.62, trim={0 0.0cm 0.25cm 0.25cm},clip]{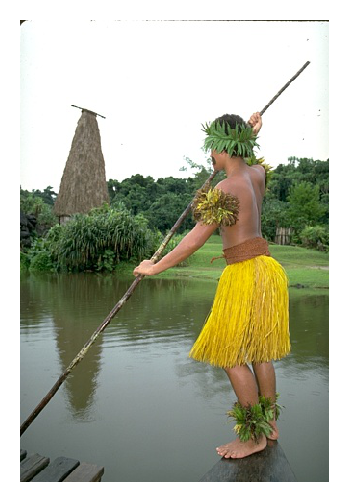}
&\includegraphics[scale=0.62, trim={0 0.0cm 0.25cm 0.25cm},clip]{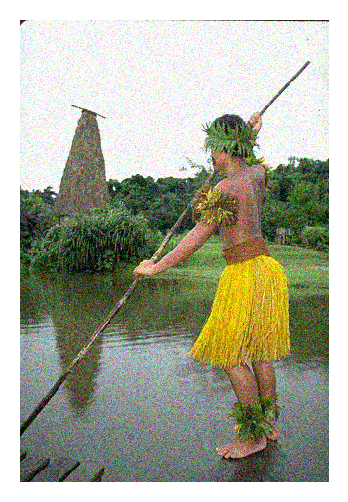}
&\includegraphics[scale=0.62, trim={0 0.0cm 0.25cm 0.25cm},clip]{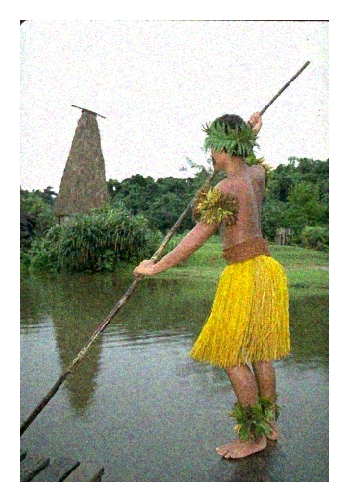}
\end{tabular}

\vspace{-0.3cm}
\caption{\small Denoising results obtained with the proposed LB-FB training strategy on two sample images from the BSDS dataset.}
    \label{fig:numerics_visual}

\vspace{-0.2cm}
\end{figure}

\begin{table}[!t]
\caption{\small Average SSIM and PSNR (and standard deviation) obtained for
the $9$ BSDS test images.}

\vspace{-0.3cm}
\begin{center}\small
\begin{tabular}{  c  c  c }
& \multicolumn{1}{c}{SSIM} & \multicolumn{1}{c}{PSNR} \\  
\hline\hline
Noisy & $0.3773 \pm 0.0929$ & $19.7534 \pm 0.3378$ \\ 
\hline
LB-FB & $0.6504 \pm 0.0874$ & $26.3502 \pm 0.4448$ \\ 
\hline
SGD & $0.6413 \pm 0.0921$ & $26.1366 \pm 0.4836$
\end{tabular}
\end{center}
\label{table}

\vspace{-0.3cm}
\end{table}

%\begin{table}[!t]
%\caption{Average SSIM and PSNR (and standard deviation) obtained for
%the test BSDS test images. \review{64x64 BS=50}}

%\vspace{-0.3cm}
%\begin{center}\small
%\begin{tabular}{  c  c  c }
%& \multicolumn{1}{c}{SSIM} & \multicolumn{1}{c}{PSNR} \\  
%\hline\hline
%Noisy & $0.3774 \pm 0.0927$ & $19.7530 \pm 0.3306$ \\ 
%\hline
%LB & $0.6576 \pm 0.0855$ & $26.5220 \pm 0.4657$ \\ 
%\hline
%SGD & $0.6496 \pm 0.0877$ & $26.3318 \pm 0.4204$
%\end{tabular}
%\end{center}
%\label{table}

%\vspace{-0.3cm}
%\end{table}

\subsection{Numerical results}

\noindent\textit{Mini-batch stability.\;}
Figure~\ref{fig:numerics_batches} gives the $\ell^2$ loss values as defined in Section~\ref{ssec:prop:LB} with respect to computational time (left) and epochs (right). %On the left hand side we plot the loss value over runtime. 
%We observe that, although all training strategies seem to converge to the same value, there is an optimal batch size in terms of runtime, equal to $10$ (or $5$).
On the one hand, we observe that using bigger batch sizes leads to slower convergence. On the other hand, using smaller batch sizes seems to enable better optimization of $\ell^2$. However the computation of the $50$ epochs takes longer runtime with very small batch sizes (e.g., $1$ or $2$). Hence a trade-off should be adopted, to avoid saturating the memory while keeping high performances.

% the full stochastic implementation takes relatively longer comparing to using the entire collection of data when the full gradient is evaluated. Using small mini-batches size such as 5 or 10 provides a balanced trade-off between per batch computation overhead and runtime. On the right we plot loss value over epochs. Under varying batch sizes, the proposed strategy shows stable performance in converging towards similar end loss values. 

\smallskip
\noindent\textit{Comparison with SGD.\;}
Figure~\ref{fig:numerics_convergence} gives the $\ell^2$ loss values with respect to epochs, obtained with the proposed LB-FB training approach and with the standard SGD, for the $3$ considered network depths. Overall, we see that LB-FB needs fewer epochs than SGD to converge, and seems more stable (no oscillation ensured by the backtracking approach). % is capable of converging in fewer to lower loss function values than SGD. We also observe that SGD approaches struggle to maintain steady decays while LB approaches exhibits more controlled decay patterns. 
We also observe that, as we increase the network depth, a light discrepancy appears between the curves. This suggests that the proposed LB-FB strategy remains stable for deeper network architectures. This is most likely due to the LB formulation, where the training objective is decomposed to layer-wise local objectives. %The increase in depth adds nothing more than additional penalty terms to the overall objective, hence allowing training to be more reactive to architectural changes. 

Finally, for visual inspection, we provide examples of denoised images obtained with our training strategy for $K=15$ in Figure~\ref{fig:numerics_visual}. %In Figure~\ref{fig:error-image}, we give error maps showing the difference between the noisy \review{punting} image and its groundtruth (left), and the difference between the denoised image obtained with \review{the proposed training strategy (resp. SGD)} and the groundtruth (middle, resp. right). 
In Table~\ref{table} we summarize PSNR and SSIM average scores and standard deviation obtained on the $9$~BSDS test images. %(middle, resp. right) $8$

\section{Conclusion}\label{sec:conclusion}

In this work, we have proposed a lifted Bregman training strategy combined with a mini-batch block-coordinate forward-backward algorithm for training an unfolded PNN specifically designed for the image Gaussian denoising task. %This strategy leads to the formulation of a bi-convex penalised loss objective and enjoys computational convenience for the evaluation of its partial gradients with respect to network parameters. To solve the learning problem, we further proposed a mini-batch stochastic proximal gradient descent algorithm. 
Through numerical experiments, we demonstrated the efficiency of the proposed strategy and showed its competitive performance in training deeper network architectures, comparing to standard SGD approaches using the back-propagation methods.

%\begin{figure}
%    \centering
%    \includegraphics[scale=0.3]{figures/Denoise_LBPNN_IMAGENET_Train_Samples_Small_NEW.png}
%    \caption{sample denoised train images with error map}
%    \label{fig:enter-label}
%\end{figure}

% \section*{Acknowledgment}

%\section*{References}

\small
\bibliographystyle{ieeetr}
\bibliography{reference}

\begin{thebibliography}{10}

\bibitem{Goodfellow-et-al-2016}
I.~Goodfellow, Y.~Bengio, and A.~Courville, {\em Deep Learning}.
\newblock MIT Press, 2016.

\bibitem{zhang2017beyond}
K.~Zhang, W.~Zuo, Y.~Chen, D.~Meng, and L.~Zhang, ``Beyond a gaussian denoiser: Residual learning of deep cnn for image denoising,'' {\em IEEE Trans. Image Process.}, vol.~26, no.~7, pp.~3142--3155, 2017.

\bibitem{ronneberger2015u}
O.~Ronneberger, P.~Fischer, and T.~Brox, ``U-net: Convolutional networks for biomedical image segmentation,'' in {\em MICCAI 2015, proceedings, part III 18}, pp.~234--241, 2015.

\bibitem{gregor2010learning}
K.~Gregor and Y.~LeCun, ``Learning fast approximations of sparse coding,'' in {\em Proc. 27th Int. Conf. Mach. Learn.}, pp.~399--406, 2010.

\bibitem{zhang2017convergence}
Z.~Zhang and M.~Brand, ``On the convergence of block coordinate descent in training dnns with tikhonov regularization,'' in {\em Adv. Neural Inf. Process. Syst.}, pp.~1719--1728, 2017.

\bibitem{hasannasab2020parseval}
M.~Hasannasab, J.~Hertrich, S.~Neumayer, G.~Plonka, S.~Setzer, and G.~Steidl, ``Parseval proximal neural networks,'' {\em J. Fourier Anal. Appl.}, vol.~26, no.~4, pp.~1--31, 2020.

\bibitem{combettes2020deep}
P.~L. Combettes and J.-C. Pesquet, ``Deep neural network structures solving variational inequalities,'' {\em Set-Valued Var. Anal.}, pp.~1--28, 2020.

\bibitem{adler2018learned}
J.~Adler and O.~{\"O}ktem, ``Learned primal-dual reconstruction,'' {\em IEEE Trans. Med. Imaging.}, vol.~37, no.~6, pp.~1322--1332, 2018.

\bibitem{jiu2021deep}
M.~Jiu and N.~Pustelnik, ``A deep primal-dual proximal network for image restoration,'' {\em IEEE J. Sel. Top. Signal Process.}, vol.~15, no.~2, pp.~190--203, 2021.

\bibitem{monga2021algorithm}
V.~Monga, Y.~Li, and Y.~C. Eldar, ``Algorithm unrolling: Interpretable, efficient deep learning for signal and image processing,'' {\em IEEE Sign. Process. Mag.}, vol.~38, no.~2, pp.~18--44, 2021.

\bibitem{kingma2014adam}
D.~P. Kingma and J.~Ba, ``Adam: A method for stochastic optimization,'' {\em arXiv preprint arXiv:1412.6980}, 2014.

\bibitem{rumelhart1986learning}
D.~E. Rumelhart, G.~E. Hinton, and R.~J. Williams, ``Learning representations by back-propagating errors,'' {\em nature}, vol.~323, no.~6088, pp.~533--536, 1986.

\bibitem{he2016deep}
K.~He, X.~Zhang, S.~Ren, and J.~Sun, ``Deep residual learning for image recognition,'' in {\em Proc. IEEE Conf. Comput. Vis. Pattern Recognit.}, pp.~770--778, 2016.

\bibitem{zach2019contrastive}
C.~Zach and V.~Estellers, ``Contrastive learning for lifted networks,'' British Machine Vision Conference (BMVC), 2019.

\bibitem{gu2020fenchel}
F.~Gu, A.~Askari, and L.~El~Ghaoui, ``Fenchel lifted networks: A lagrange relaxation of neural network training,'' in {\em Int. Conf. Artif. Intell. Stat.}, pp.~3362--3371, PMLR, 2020.

\bibitem{taylor2016training}
G.~Taylor, R.~Burmeister, Z.~Xu, B.~Singh, A.~Patel, and T.~Goldstein, ``Training neural networks without gradients: A scalable admm approach,'' in {\em Int. Conf. Mach. Learn.}, pp.~2722--2731, 2016.

\bibitem{carreira2014distributed}
M.~Carreira-Perpinan and W.~Wang, ``Distributed optimization of deeply nested systems,'' in {\em Artificial Intelligence and Statistics}, pp.~10--19, PMLR, 2014.

\bibitem{xu2022alternative}
C.~Xu, X.~Cheng, and Y.~Xie, ``An alternative approach to train neural networks using monotone variational inequality,'' {\em arXiv preprint arXiv:2202.08876}, 2022.

\bibitem{frecon2022bregman}
J.~Frecon, G.~Gasso, M.~Pontil, and S.~Salzo, ``Bregman neural networks,'' in {\em Int. Conf. Mach. Learn.}, pp.~6779--6792, 2022.

\bibitem{combettes2023variational}
P.~L. Combettes, J.-C. Pesquet, and A.~Repetti, ``A variational inequality model for learning neural networks,'' in {\em Proc. IEEE Int. Conf. Acous. Speech Sig. Process. (ICASSP) 2023}, pp.~1--5, 2023.

\bibitem{wang2023lifted}
X.~Wang and M.~Benning, ``Lifted bregman training of neural networks,'' {\em J. Mach. Learn. Res.}, vol.~24, no.~232, pp.~1--51, 2023.

\bibitem{bregman1967relaxation}
L.~M. Bregman, ``The relaxation method of finding the common point of convex sets and its application to the solution of problems in convex programming,'' {\em USSR Comput. Math. Math. Phys.}, vol.~7, no.~3, pp.~200--217, 1967.

\bibitem{kiwiel1997proximal}
K.~C. Kiwiel, ``Proximal minimization methods with generalized bregman functions,'' {\em SIAM J. Control Optim.}, vol.~35, no.~4, pp.~1142--1168, 1997.

\bibitem{rudin1992nonlinear}
L.~I. Rudin, S.~Osher, and E.~Fatemi, ``Nonlinear total variation based noise removal algorithms,'' {\em Physica D Nonlinear Phenom.}, vol.~60, no.~1-4, pp.~259--268, 1992.

\bibitem{mallat1999wavelet}
S.~Mallat, {\em A wavelet tour of signal processing}.
\newblock Elsevier, 1999.

\bibitem{jacques2011panorama}
L.~Jacques, L.~Duval, C.~Chaux, and G.~Peyr{\'e}, ``A panorama on multiscale geometric representations, intertwining spatial, directional and frequency selectivity,'' {\em Signal Process.}, vol.~91, no.~12, pp.~2699--2730, 2011.

\bibitem{combettes2010dualization}
P.~L. Combettes, {\DJ}.~D{\~u}ng, and B.~C. V{\~u}, ``Dualization of signal recovery problems,'' {\em Set-Valued Var. Anal.}, vol.~18, no.~3, pp.~373--404, 2010.

\bibitem{repetti2022dual}
A.~Repetti, M.~Terris, Y.~Wiaux, and J.-C. Pesquet, ``Dual forward-backward unfolded network for flexible plug-and-play,'' in {\em 30th Eur. Signal Process. Conf. (EUSIPCO)}, pp.~957--961, 2022.

\bibitem{le2022faster}
H.~T.~V. Le, N.~Pustelnik, and M.~Foare, ``The faster proximal algorithm, the better unfolded deep learning architecture? the study case of image denoising,'' in {\em 30th Eur. Signal Process. Conf. (EUSIPCO)}, pp.~947--951, 2022.

\bibitem{le2023pnn}
H.~T.~V. Le, A.~Repetti, and N.~Pustelnik, ``Pnn: From proximal algorithms to robust unfolded image denoising networks and plug-and-play methods,'' {\em arXiv preprint arXiv:2308.03139}, 2023.

\bibitem{zhang2021plug}
K.~Zhang, Y.~Li, W.~Zuo, L.~Zhang, L.~Van~Gool, and R.~Timofte, ``Plug-and-play image restoration with deep denoiser prior,'' {\em IEEE Transactions on Pattern Analysis and Machine Intelligence}, vol.~44, no.~10, pp.~6360--6376, 2021.

\bibitem{lecun1988theoretical}
Y.~LeCun, D.~Touresky, G.~Hinton, and T.~Sejnowski, ``A theoretical framework for back-propagation,'' in {\em Proc. 1988 Connect. Models Summer Sch.}, vol.~1, pp.~21--28, 1988.

\bibitem{wang2023inversion}
X.~Wang and M.~Benning, ``A lifted bregman formulation for the inversion of deep neural networks,'' {\em Front. Appl. Math. Stat.}, vol.~9, p.~1176850, 2023.

\bibitem{chouzenoux2016block}
E.~Chouzenoux, J.-C. Pesquet, and A.~Repetti, ``A block coordinate variable metric forward--backward algorithm,'' {\em J. Glob. Optim.}, vol.~66, no.~3, pp.~457--485, 2016.

\bibitem{russakovsky2015imagenet}
O.~Russakovsky, J.~Deng, H.~Su, J.~Krause, S.~Satheesh, S.~Ma, Z.~Huang, A.~Karpathy, A.~Khosla, M.~Bernstein, {\em et~al.}, ``Imagenet large scale visual recognition challenge,'' {\em Int. J. Comput. Vis.}, vol.~115, pp.~211--252, 2015.

\bibitem{martin2001database}
D.~Martin, C.~Fowlkes, D.~Tal, and J.~Malik, ``A database of human segmented natural images and its application to evaluating segmentation algorithms and measuring ecological statistics,'' in {\em Proc. 8th IEEE Int. Conf. Comput. Vis.}, vol.~2, pp.~416--423, 2001.

\bibitem{glorot2010understanding}
X.~Glorot and Y.~Bengio, ``Understanding the difficulty of training deep feedforward neural networks,'' in {\em Proc. 13th Int. Conf. Artif. Intell. Stat.}, pp.~249--256, JMLR Workshop and Conference Proceedings, 2010.

\end{thebibliography}

%\begin{thebibliography}{00}
%\bibitem{b1} G. Eason, B. Noble, and I. N. Sneddon, ``On certain integrals of Lipschitz-Hankel type involving products of Bessel functions,'' Phil. Trans. Roy. Soc. London, vol. A247, pp. 529--551, April 1955.
%\bibitem{b2} J. Clerk Maxwell, A Treatise on Electricity and Magnetism, 3rd ed., vol. 2. Oxford: Clarendon, 1892, pp.68--73.
%\bibitem{b3} I. S. Jacobs and C. P. Bean, ``Fine particles, thin films and exchange anisotropy,'' in Magnetism, vol. III, G. T. Rado and H. Suhl, Eds. New York: Academic, 1963, pp. 271--350.

%\end{thebibliography}

\end{document}